\documentclass[a4paper,11pt]{amsart}
\usepackage{amsmath,amssymb,amsfonts}

\theoremstyle{theorem}
\newtheorem{theorem}{Theorem}[section]
\newtheorem{lemma}{Lemma}[section]
\newtheorem{corollary}{Corollary}[section]
\newtheorem{definition}{Definition}[section]
\theoremstyle{definition}
\newtheorem{example}{Example}[section]

\begin{document}

\title[New properties for a composition of some generating functions]{New properties for a composition of some generating functions for primes}
\author{Dmitry V. Kruchinin, Yuriy V. Shablya}
\address{Tomsk State University of Control Systems and Radioelectronics, Tomsk, Russia}
\email{kruchininDm@gmail.com}


\begin{abstract}
In this paper, we consider properties of coefficients of a generating functions composition, where the outer function is a logarithmic generating function and the inner function is an ordinary generating function with integer coefficients.
Using notions of composita and composition of generating functions, we get new properties for this composition.
The properties can be used for distinguishing prime numbers from composite numbers.
As an application, obtained results can be used to obtain new primality criteria.
We obtain primality criteria for the Mersenne numbers, the Lucas numbers, the Pell-Lucas numbers, the Jacobsthal-Lucas numbers, and the Lucas sequences.

\vspace{2mm}


\end{abstract}


\maketitle


\section{INTRODUCTION}

There are many authors who have studied generating functions and their properties (for instance, Comtet \cite{Comtet1974}, Flajolet \cite{Flajolet2009}, Graham \cite{Graham1989}, Robert \cite{Robert2000}, Stanley \cite{Stanley1978}, Wilf \cite{Wilf1994}).
Generating functions are a powerful tool for solving problems in number theory, combinatorics, algebra, probability theory, and in other fields of mathematics.
One of the advantages of generating functions is that an infinite number sequence can be represented in the form of single expression.

\begin{definition}
     The ordinary generating function of a sequence $\{a_n\}$ is the following formal power series
     $$
          A(x) = a_0 + a_1 x + a_2 x^2 + \ldots = \sum_{n \geq 0} {a_n x^n}.
     $$
\end{definition}

There exist many operations on generating functions, such as addition, multiplication, integration, differentiation, and composition.
In this paper, we consider properties of a composition of generating functions, where the outer function is a logarithmic generating function and the inner function is an ordinary generating function with integer coefficients.

The article structure is presented as follows:
\begin{itemize}
\item In section~2, we show basic rules for calculating the coefficients of the composition of generating functions and we give necessary mathematical notions.
\item In section~3 we get the main results and give them in Theorem~\ref{theorem1}, Theorem~\ref{theorem2} and Corollary~\ref{corollary1}.
\item In section~4 we consider an application of obtained results.
\end{itemize}

\section{COMPOSITION OF GENERATING FUNCTIONS}

According to Stanley \cite{Stanley1999}, the composition of generating functions is defined as follows:

\begin{definition}
     Suppose $R(x) = \sum_{n \geq 0} {r_n x^n}$ and $F(x) = \sum_{n > 0} {f_n x^n}$ are formal power series. Then the composition $G(x) = R\left(F(x)\right) = \sum_{n \geq 0} {r_n F(x)^n}$ is a well-defined formal power series.
\end{definition}

We consider a composition of generating functions $R\left(F(x)\right)$, where the outer function $R(x)$ is a logarithmic generating function and the inner function $F(x)$ is an ordinary generating function with integer coefficients.

According to Robert \cite{Robert2000}, the logarithmic generating function is defined as follows:
\begin{definition}
     The logarithmic generating function is the generating function in the form
     $$
          R(x) = \sum_{n > 0} {r_n x^n} = \sum_{n > 0} {\frac{a_n}{n} x^n}.
     $$
\end{definition}

To calculate the coefficients of the composition of generating functions, we use the mathematical notion of the composita of a given generating function, which was introduced by Kruchinin \cite{Kruchinin2014}.

\begin{definition}
     The composita of the generating function $F(x) = \sum_{n > 0} {f_n x^n}$ is the function of two variables
     \begin{equation}
          \label{Composita}
          F^{\Delta}(n, k) = \sum_{\pi_k \in C_n} {f_{\lambda_1} f_{\lambda_2} \ldots f_{\lambda_k}},
     \end{equation}
     where $C_n$ is a set of all compositions of an integer $n$, $\pi_k$ is the composition $n$ into $k$ parts such that $\sum_{i = 1}^{k} {\lambda_i} = n$.
\end{definition}

Also we can write the following condition
$$
     \left(F(x)\right)^k = \sum_{n \geq k} {  F^{\Delta}(n, k) x^n }.
$$

Using the notion of the composita \eqref{Composita}, we can obtain the coefficients of the composition of generating functions \cite{Kruchinin2014}.

\begin{lemma}
     Suppose $F(x) = \sum_{n > 0} {f_n x^n}$ and $R(x) = \sum_{n \geq 0} {r_n x^n}$ are generating functions, and $F^{\Delta}(n, k)$ is the composita of $F(x)$.
     Then for the composition of the generating functions $G(x) = R\left(F(x)\right)$ the coefficients of the generating function $G(x) = \sum_{n \geq 0} {g_n x^n}$ are
     \begin{equation}
          \label{CompFormula}
          g_n = 
          \begin{cases}
               r_0, & \text{if $n = 0$;} \\
               \sum_{k = 1}^{n} {F^{\Delta}(n, k) r_k}, & \text{if $n > 0$.}
          \end{cases}
     \end{equation}
\end{lemma}

\section{MAIN RESULTS}

We show the main results of this paper in the following theorems and their corollaries.
\begin{theorem}
     \label{theorem1}
     Suppose $F(x) = \sum_{n > 0} {f_n x^n}$ is an ordinary generating function, $R(x) = \sum_{n > 0} {r_n x^n} = \sum_{n > 0} {\frac{a_n}{n} x^n}$ is a logarithmic generating function, $\{f_n\}$ and $\{a_n\}$ are integer sequences, $G(x) = \sum_{n > 0} {g_n x^n}$ is a generating function, which is the composition of the generating functions $G(x) = R\left(F(x)\right)$.
     Then the value of
     \begin{equation}
          \label{theorem1_ngn}
          n g_n = n \sum_{k = 1}^{n} {F^{\Delta}(n, k)  r_k} = n \sum_{k = 1}^{n} {\frac{F^{\Delta}(n, k) a_k}{k}}
     \end{equation}
     is integer for $n \in \mathbb{N}$.
\end{theorem}

\begin{proof}
     According to \eqref{CompFormula}, the coefficients of the generating function $G(x)$ are 
     \begin{equation}
          \label{gn}
          g_n = \sum_{k = 1}^{n} {F^{\Delta}(n, k)  r_k} = \sum_{k = 1}^{n} {\frac{F^{\Delta}(n, k) a_k}{k}}.
     \end{equation}

     The derivative of the generating function $G(x) = \sum_{n > 0} {g_n x^n}$ (see \cite{Lando2003}) is the function 
     $$
          G'(x) = g_1 + 2 g_2 x^2 + 3 g_3 x^2 + \ldots + n g_n x^{n - 1} + \ldots = \sum_{n > 0} {n g_n x^{n - 1}}.
     $$
     
     Also, we can obtain the derivative of the generating function $G(x)$
     $$
          G'(x) = \left(R\left(F(x)\right)\right)' = F'(x) R'\left(F(x)\right).
     $$
     
     Since $\{n f_n\}$ is an integer sequence, then the coefficients of the generating function $F'(x) = \sum_{n > 0} {n f_n x^{n - 1}}$ are integer.
     
     Since $\{a_n\}$ is an integer sequence, then the coefficients of the generating function $R'(x) = \sum_{n > 0} {n r_n x^{n - 1}} = \sum_{n > 0} {a_n x^{n - 1}}$ are integer.
     
     Since $F(x)$ and $R'(x)$ are generating functions with integer coefficients, then the coefficients of the composition of the generating functions $R'\left(F(x)\right)$ are also integer.
     
     The coefficients of the product of generating functions with integer coefficients are also integer, so $G'(x) = F'(x) R'\left(F(x)\right) = \sum_{n > 0} {n g_n x^{n - 1}}$ is the generating function with integer coefficients.
     
     Therefore, $\{n g_n\}$ is an integer sequence.
     
     The theorem is proved: the value of the expression \eqref{theorem1_ngn} is integer.
\end{proof}

From Theorem~\ref{theorem1} we obtain a new important property of the coefficients of the generating functions composition:

\begin{corollary}
     \label{corollary1}
     The value of the coefficients function \eqref{gn} without the $n$-th term
     \begin{equation}
          \label{Criterion}
           \sum_{k = 1}^{n - 1} {\frac{F^{\Delta}(n, k) a_k}{k}} = \frac{n g_n - a_n f_1^n}{n}
     \end{equation}
     is integer for every prime $n$.
     The converse is not true.
\end{corollary}
 	
\begin{proof}
     Let us consider the expression \eqref{theorem1_ngn}.
     The value of the $n$-th term in the sum \eqref{theorem1_ngn} is equal to the integer expression
     $$
          n \frac{F^{\Delta}(n, n) a_n}{n} = a_n f_1^n.
     $$
     
     Then the value of the expression
     $$
          n g_n - a_n f_1^n = n \sum_{k = 1}^{n - 1} {\frac{F^{\Delta}(n, k) a_k}{k}}
     $$
     is integer.
      
     The coefficients of the generating function $F(x) = \sum_{n > 0} {f_n x^{n - 1}}$ are integer.
     Then the values of the composita $F^{\Delta}(n, k)$ are integer.
     
     Since $n$ is a prime, $n > k$, $F^{\Delta}(n, k)$ and $a_k$ are integer, then the expression \eqref{Criterion} is integer.
\end{proof}

The expression \eqref{Criterion} is integer for every prime $n$.
It can be used for distinguishing prime numbers from composite numbers.
For composite numbers $n$ the value of the expression \eqref{Criterion} can be integer or not integer.
But if the value of the expression \eqref{Criterion} is not integer, then $n$ is determinately a composite number.

We note that an integral of the generating function $B(x) = \sum_{n \geq 0}{b_n x^n}$ through term by term integration of power series (see \cite{Lando2003}) is the function 
\begin{equation}
     \label{intBogf}
     \int{B(x)} = b_0 x + b_1 \frac{x^2}{2} + b_2 \frac{x^3}{3} + \ldots + b_n \frac{x^{n + 1}}{n + 1} + \ldots = \sum_{n > 0} {\frac{b_{n - 1}}{n} x^n}.
\end{equation}

Then Theorem~\ref{theorem1} takes the following form:

\begin{theorem}
     \label{theorem2}
     Suppose $B(x) = \sum_{n \geq 0}{b_n x^n}$ and $F(x) = \sum_{n > 0} {f_n x^n}$ are ordinary generating functions with integer coefficients, $G(x) = \sum_{n > 0} {g_n x^n}$ is a generating function, that is the composition of the generating functions $G(x) = R\left(F(x)\right)$, where $R(x) = \int{B(x)}$.
     Then the value of
     $$
          n g_n = n \sum_{k = 1}^{n} {\frac{F^{\Delta}(n, k) b_{k - 1}}{k}}
     $$
     is integer for $n \in \mathbb{N}$.
\end{theorem}

\begin{proof}
     If in \eqref{intBogf} we consider a sequence $\{a_n\}$, where $a_n = b_{n - 1}$ for every integer $n = 1, 2, \ldots$, then we obtain the logarithmic generating function
     $$
          \int{B(x)} = \sum_{n > 0} {\frac{b_{n - 1}}{n} x^n} = \sum_{n > 0} {\frac{a_n}{n} x^n} = \sum_{n > 0} {r_n x^n} = R(x).
     $$
     
     According to Theorem~\ref{theorem1}, the value of the expression
     $$
          n g_n = n \sum_{k = 1}^{n} {\frac{F^{\Delta}(n, k) a_k}{k}} = n \sum_{k = 1}^{n} {\frac{F^{\Delta}(n, k) b_{k - 1}}{k}}
     $$
     is integer for $n \in \mathbb{N}$.
\end{proof}

     Also, the value of the expression
     $$
           \sum_{k = 1}^{n - 1} {\frac{F^{\Delta}(n, k) b_{k - 1}}{k}} = \frac{n g_n - b_{n - 1} f_1^n}{n}
     $$
     is integer for every prime $n$.
     The converse is not true.

\section{APPLICATION}

The obtained results can be used in public-key cryptography \cite{Song2009}.
The obtained properties of the composition of generating functions can be used in constructing new primality criteria.
The primality criterion is a statement that necessarily must be satisfied for primes.
Such criteria may be used as the basis for the probabilistic primality tests.
If the number does not satisfy the conditions of tests, then probabilistic primality tests can exactly detect that the number $n$ is composite.
If the number satisfies all conditions, then $n$ is a prime with a certain probability.
Often repeating the test with different parameters can reduce the probability of error.
Using the obtained properties of the composition of generating functions, we consider several examples of construction of primality criteria.

\begin{example}
     Let us consider the following composition of the generating functions $G(x) = R\left(F(x)\right)$, where
     \begin{equation}
          \label{Ex1_Fx}
          F(x) = \frac{b x}{1 - a x} = \sum_{n > 0} {b a^{n - 1} x^n}
     \end{equation}
     is a generating function with integer coefficients, $b$ and $a$ are integer;
     $$
          R(x) = \ln\left(\frac{1}{1 - x}\right) = \sum_{n > 0} {\frac{1}{n} x^n}
     $$
     is a logarithmic generating function. 

     Since the composita of the generating function \eqref{Ex1_Fx} is (see \cite{Kruchinin2013})
     $$
          F^{\Delta}(n, k, a, b) = \binom{n - 1}{k - 1} a^{n - k} b^k,
     $$
     then the coefficients of the composition $G(x) = R\left(F(x)\right)$ are given by
     \begin{equation}
          \label{Ex1_gn}
          g_n = \sum_{k = 1}^{n} {\binom{n - 1}{k - 1} \frac{a^{n - k} b^k}{k}}.
     \end{equation}

     Using \eqref{Criterion}, we get the primality criterion: the value of expression
     $$
          \sum_{k = 1}^{n - 1} {\binom{n - 1}{k - 1} \frac{a^{n - k} b^k}{k}} = \frac{(a + b)^n - a^n - b^n}{n}.
     $$
     is integer for every prime $n$.
     
     After some transformations, we obtain the following primality criterion: if $n$ is prime, then
     \begin{equation}
          \label{Ex1_Criterion2}
          (a + b)^n - a^n - b^n \equiv 0 \mod n.
     \end{equation}
     
     Obtained expression is similar to the Fermat's little theorem \cite{Agrawal2006}.
     
     Since the following condition
     $$
          a^n \equiv a \mod n
     $$
     holds true for every prime $n$ and integer $a$, then we transform \eqref{Ex1_Criterion2} into the following primality criterion: if $n$ is prime, then
     \begin{equation}
          \label{Ex1_Criterion3}
          (a + b)^n \equiv a + b \mod n.
     \end{equation}

     Substituting $c$ for $(a + b)$ in \eqref{Ex1_Criterion3}, we get the Fermat's little theorem: if $n$ is prime, then for every $c \in \{1, 2, \ldots, n - 1\}$
     $$
          c^{n - 1} \equiv 1 \mod n.
     $$

     The Fermat's little theorem has had a great influence in algorithmic number theory as it has been the basis for some of the most well-known primality tests: Fermat primality test, Solovay-Strassen primality test, Miller-Rabin primality test, AKS primality test, etc. (see \cite{Agrawal2006}).
     The first three are probabilistic polynomial time algorithms and are widely used in practice, the fourth one is the only known deterministic polynomial time algorithm.
     
     If we set $a = 1$ and $b = 1$ in \eqref{Ex1_gn}, then we obtain
     $$
          n g_n = \{1, 3, 7, 15, 31, 63, 127, 255, 511, 1023, \ldots\} = 2^n - 1 = M_n,
     $$
     where $M_n$ is the Mersenne number (the sequence A000225 in The On-Line Encyclopedia of Integer Sequences \cite{OEIS}).
     
     Therefore, using Corollary~\ref{corollary1}, we obtain the primality criterion: if $n$ is prime, then
     $$
          M_n \equiv 1 \mod n.
     $$
\end{example}

\begin{example}
     Let us consider the following composition of generating functions $G(x) = R\left(F(x)\right)$, where
     \begin{equation}
          \label{Ex2_Fx}
          F(x) = a x + b x^2
     \end{equation}
     is a generating function with integer coefficients, $b$ and $a$ are integer;
     $$
          R(x) = \ln\left(\frac{1}{1 - x}\right) = \sum_{n > 0} {\frac{1}{n} x^n}
     $$
     is a logarithmic generating function. 

     Since the composita of the generating function \eqref{Ex2_Fx} is (see \cite{Kruchinin2013})
     $$
          F^{\Delta}(n, k, a, b) = \binom{k}{n - k} a^{2 k - n} b^{n - k},
     $$
     then the coefficients of the composition $G(x) = R\left(F(x)\right)$ are given by
     \begin{equation}
          \label{Ex2_gn}
          g_n =\sum_{k = 1}^{n} {\binom{k}{n - k} \frac{a^{2 k - n} b^{n - k}}{k}}.
     \end{equation}

     Using \eqref{Criterion}, we get the new primality criterion: the value of expression
     \begin{equation}
          \label{Ex2_Criterion}
          \sum_{k = 1}^{n - 1} {\binom{n - 1}{k - 1} \frac{a^{n - k} b^k}{k}} = \frac{\left(\frac{a + \sqrt{a^2 + 4 b}}{2}\right)^n + \left(\frac{a - \sqrt{a^2 + 4 b}}{2}\right)^n - a^n}{n}
     \end{equation}
     is integer for every prime $n$.
     
     We consider special cases of $a$ and $b$.
     
     Set $a = 1$ and $b = 1$ in \eqref{Ex2_gn}, then we obtain
     $$
          n g_n = \{1, 3, 4, 7, 11, 18, 29, 47, 76, \ldots\} = \left(\frac{1 + \sqrt{5}}{2}\right)^n + \left(\frac{1 - \sqrt{5}}{2}\right)^n = L_n,
     $$
     where $L_n$ is the Lucas number (the sequence A000032 in The On-Line Encyclopedia of Integer Sequences \cite{OEIS}).
     
     Therefore, the primality criterion: if $n$ is prime, then
     $$
          L_n \equiv 1 \mod n.
     $$
     
     Set $a = 2$ and $b = 1$ in \eqref{Ex2_gn}, then we obtain
     $$
          n g_n = \{2, 6, 14, 34, 82, 198, 478, 1154, 2786, \ldots\} = (1 + \sqrt{2})^n + (1 - \sqrt{2})^n = Q_n,
     $$
     where $Q_n$ is the Pell-Lucas number (the sequence A002203 in The On-Line Encyclopedia of Integer Sequences \cite{OEIS}).
     
     Therefore, the primality criterion: if $n$ is prime, then
     $$
          Q_n \equiv 2 \mod n.
     $$
     
     Set $a = 1$ and $b = 2$ in \eqref{Ex2_gn}, then we obtain
     $$
          n g_n = \{1, 5, 7, 17, 31, 65, 127, 257, 511, 1025, \ldots\} = 2^n + (-1)^n = j_n,
     $$
     where $j_n$ is the Jacobsthal-Lucas number (the sequence A014551 in The On-Line Encyclopedia of Integer Sequences \cite{OEIS}).
     
     Therefore, the primality criterion: if $n$ is prime, then
     $$
          j_n \equiv 1 \mod n.
     $$
     
     Next we consider the Lucas sequence \cite{Ribenboiml2004}
     $$
          V_n(P, Q) = \left(\frac{P + \sqrt{P^2 - 4 Q}}{2}\right)^n + \left(\frac{P - \sqrt{P^2 - 4 Q}}{2}\right)^n,
     $$
     and set $a = P$, $b = -Q$ in expression \eqref{Ex2_Criterion}, then we obtain the primality criterion: the value of expression
     $$
          \frac{V_n(P, Q) - P^n}{n}
     $$
     is integer for every prime $n$.
     
     After some transformations, we obtain the following primality criterion: if $n$ is prime, then
     $$
          V_n(P, Q) \equiv P \mod n.
     $$
\end{example}

\end{document}